\documentclass[twoside,12pt]{article}
\usepackage{graphics}
\usepackage{graphicx}
\usepackage{amssymb,amsfonts}
\usepackage{amsmath}

\title{\bf A geometrical proof of sum of $\cos n\varphi$ 
\thanks{The final publication is available at {\sc N\'{e}meth, L.}, {A geometrical proof of sum of $\cos n\varphi$}, {\it Studies of the University of \v{Z}ilina}, Mathematical Series,  \textbf{27} (2015) 63-66.}
}
\author{L\'aszl\'o N\'emeth}
\date{}

\begin{document}

\maketitle

\begin{abstract}
In this article, we present a geometrical proof of sum of $\cos \ell\varphi$ where $\ell$ goes from $1$ up to $m$. Although there exist some summation forms and the proofs are simple,  they use complex numbers. Our proof comes from a geometrical construction. Moreover, from this geometrical construction we obtain an other summation form. \\[1mm]

\end{abstract}

\section{Introduction}\label{sec:introduction}

The Lagrange's trigonometric identities  are well-known formulas. The one for sum of $\cos \ell\varphi$ $(\varphi\in (0,2\pi))$ is
\begin{eqnarray}\label{eq:sum}
\sum_{\ell=1}^{m}\cos \ell\varphi=\frac12\left(\frac{\sin\,(m+\frac12)\varphi}{\sin\frac12\varphi}-1 \right). 
\end{eqnarray}
The proof of equation \eqref{eq:sum} is based on the theorem of the complex numbers in all the books, articles and lessons at the universities (\cite{Muniz}, \cite{Jeff}). In the following we give a geometrical construction which implies the formula \eqref{eq:sum} and using certain geometrical properties we obtain an other summation formula without half angles $(\varphi\in (0,2\pi),\varphi\neq \pi)$
\begin{eqnarray}\label{eq:sum2}
\sum_{\ell=1}^{m}\cos \ell\varphi  
&=&  \frac12 \left( \frac{\sin (m+1)\varphi+\sin m\varphi}{\sin\varphi}  -1 \right).
\end{eqnarray}

\section{Geometrical construction}\label{sec:geom}

\begin{figure}[!htb]
 \centering   \includegraphics{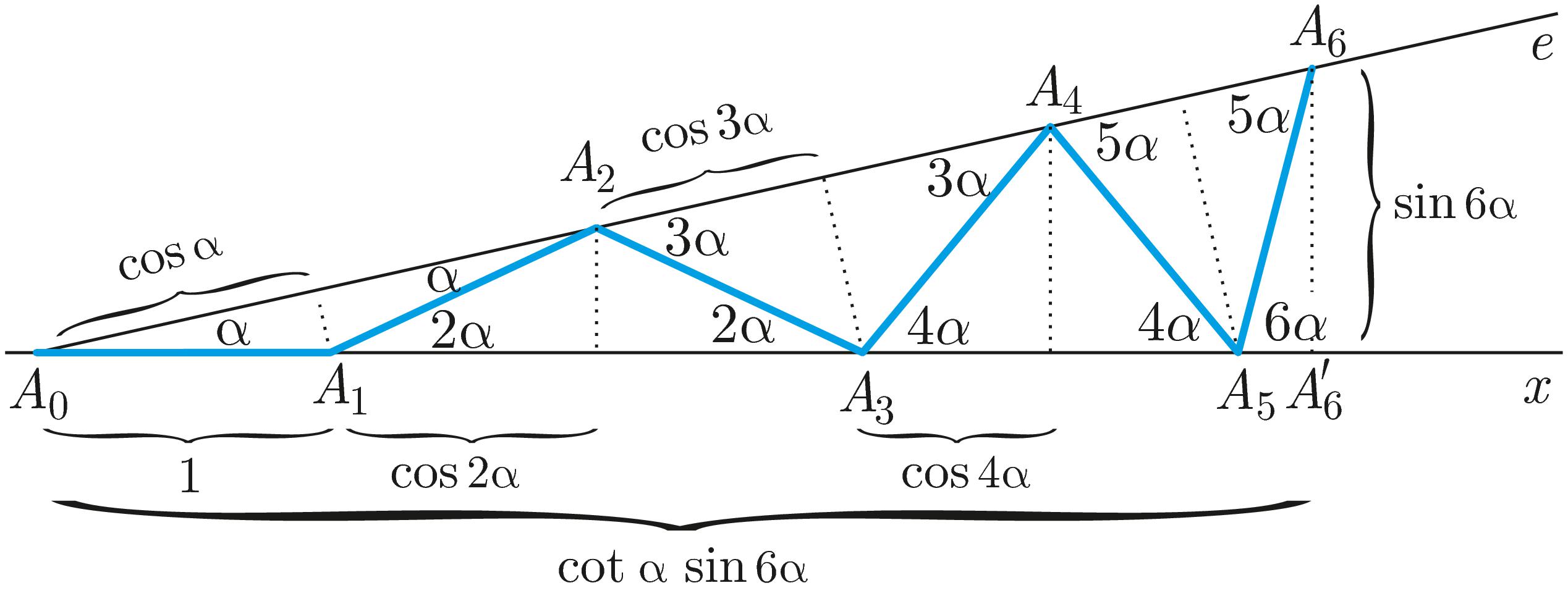}
 \caption{First seven points of the geometrical construction} \label{fig:even}
\end{figure}

\begin{figure}[!htb]
\centering\includegraphics{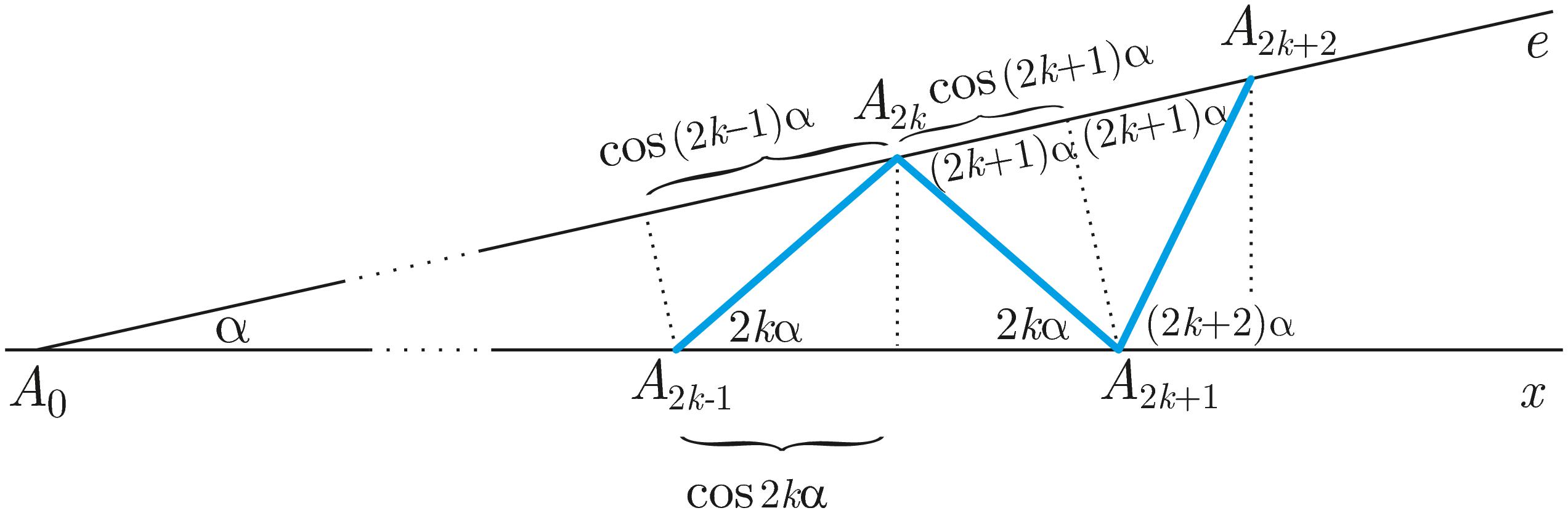}
\caption{General points of the geometrical construction}
\label{fig:even_n}
\end{figure}

\begin{figure}[!htb]
\centering\includegraphics{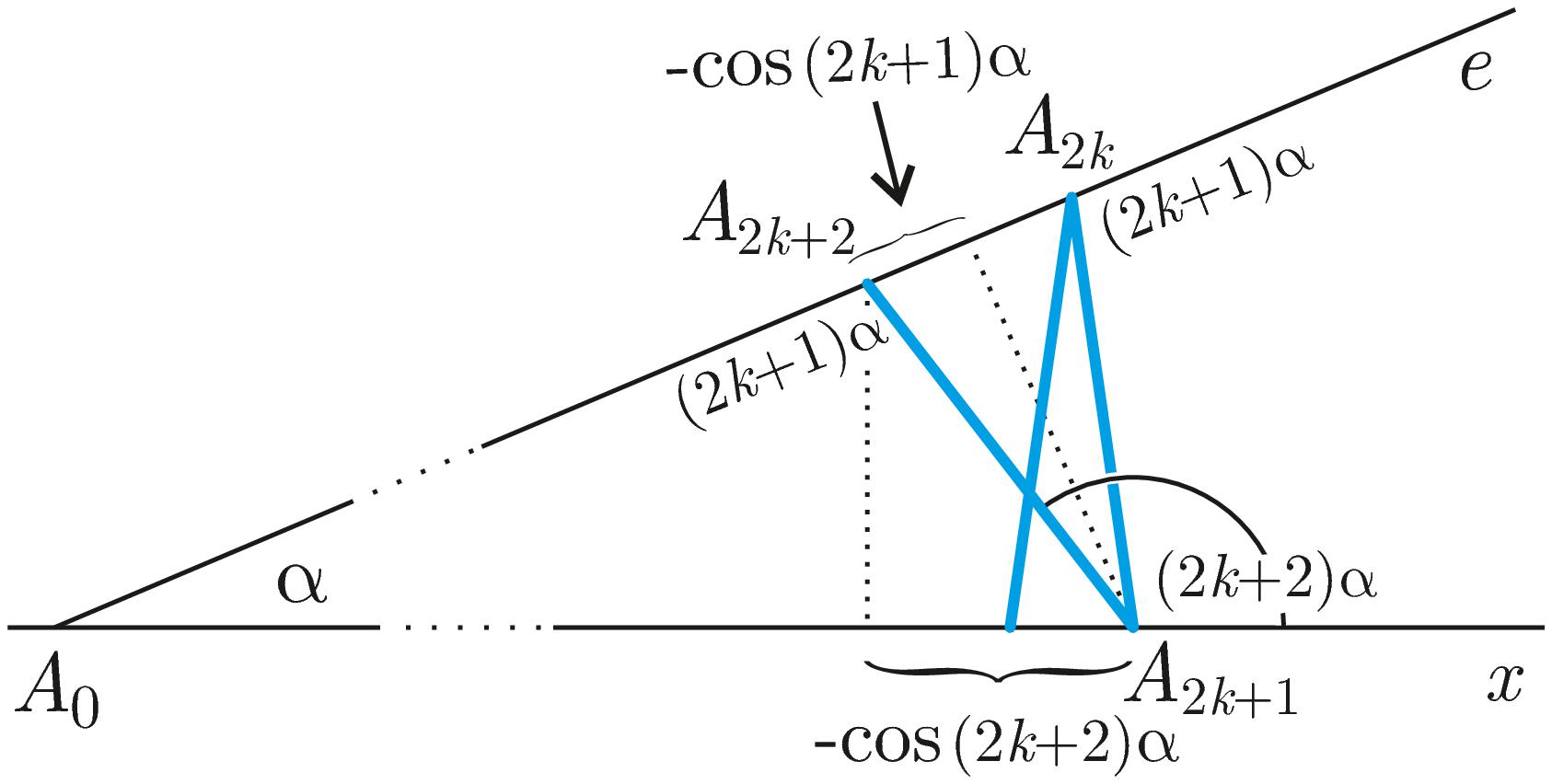}
\caption{General points}
\label{fig:even_n2}
\end{figure}

Let $x$ and $e$ be two lines with the intersection point $A_0$. Let the angle of them is $\alpha$ as $x$ is rotated to  $e$ (Figure~\ref{fig:even}). 
Let the point $A_1$ be given on $x$ such that the distance between the points $A_0$ and $A_1$ is $1$. Let the point $A_2$ be on the line $e$ such that the distance of $A_1$ and $A_2$ is also equal to $1$ and $A_2\neq A_0$ if $\alpha\ne \pi/2$ and  $\alpha\ne 3\pi/2$. Then let the new point $A_3$ be on the line $x$ again such that $A_2 A_3=1$ and $A_3\neq A_1$ if it is possible. Recursively, we can define the point $A_{\ell}$ $(\ell\geq2)$ on one of the lines $x$ or $e$ if ${\ell}$ is odd or even, respectively, where $A_{\ell-1} A_{\ell}=1$ and $A_{\ell}\neq A_{\ell-2}$ if it is possible.   Figure~\ref{fig:even} shows the first six points and Figure~\ref{fig:even_n} shows some general points. 
We can easily check that the rotation angels at vertices $A_{\ell}$ $(\ell\geq1)$ between the line $e$ (or the axis $x$) and the segments $A_{{\ell}-1}A_{\ell}$ or $A_{\ell} A_{\ell+1}$ are $(\ell-1)\alpha$ or $(\ell+1)\alpha$, respectively, as the triangles $A_{\ell-1}A_{\ell}A_{\ell+1}$ are isosceles.
(The angle  $i\alpha$ can be larger the $\pi/2$, even larger than $2\pi$. The vertices $A_\ell$ can be closer to $A_0$ then $A_{\ell-2}$ -- see Figure~\ref{fig:even_n2}.)
If $A_1$ is on the line $e$ we obtain a similar geometric construction. In that case those points $A_{\ell}$ are on line $x$ which have even indexes.   

Let $A_0$ be the origin and the line $x$ is the axis $x$.  Then the equation of the line $e$ is $\cos\alpha \cdot y=\sin \alpha \cdot x$. Let $A'_n$ be the orthogonal projection of $A_n\in e$ $(n\geq 2)$ onto the axis $x$ then from right angle triangle $A_0A'_nA_n$  the coordinates of the points $A_n$ (see Figure~\ref{fig:even}) are 
\begin{equation}\label{eq:lem00}
\begin{aligned}
  x_{n}(\alpha)&= \cot\alpha \sin n\alpha,\\
  y_{n}(\alpha)&= \sin n\alpha.
\end{aligned}
\end{equation}

If $\alpha=\pi/2$ and $\alpha=3\pi/2$ then all the points $A_n$ coincide the points $A_0$ or $A_1$, so in the following we exclude this cases.   
  
The parametric equation system of the orbits of the points $A_n\in e$ can be given  by the help of the Chebyshev polynomial too (for more details and for some figures of orbits, see in \cite{nem-lem}). The equation system is
\begin{equation}\label{eq:lem01}
\begin{aligned}
  x_{n}(\alpha)&=& \cos \alpha\ U_{n-1}( \cos \alpha), \\
  y_{n}(\alpha)&=& \sin \alpha\ U_{n-1}( \cos \alpha),
\end{aligned}
\end{equation}
where $\alpha$ goes from $0$ to $2\pi$
and $U_{n-1}(x)$ is a Chebyshev polynomial of the second kind \cite{riv}.  

Let $A_1\in x$ and $n$ be even so that $n=2k+2$.   We take the orthogonal projections of the segments $A_{{\ell}-1}A_{\ell}$ $({\ell}=1,2,...,n)$ onto the line $x$ (see Figure \ref{fig:even} and \ref{fig:even_n}). Then we realize 
\begin{eqnarray}\label{eq:even01}
x_{n}(\alpha)&=&1+2\left( \cos 2\alpha + \cos 4\alpha+ \cdots +\cos 2k\alpha \right)+ \cos (2k+2)\alpha ,
\end{eqnarray}
on the other hand, from \eqref{eq:lem00} or from  \eqref{eq:lem01} we have
\begin{eqnarray}\label{eq:even02}
x_{n}(\alpha) = \cos\alpha\ \frac{\sin(2k+2)\alpha}{\sin\alpha} .
\end{eqnarray}

Comparing \eqref{eq:even01} and \eqref{eq:even02} we obtain
\begin{eqnarray}\label{eq:even03}
1+2\sum_{{\ell}=1}^{k}\cos 2{\ell}\alpha + \cos (2k+2)\alpha
&=& \cos\alpha\ \frac{\sin(2k+2)\alpha}{\sin\alpha} .
\end{eqnarray}
 
Using the addition formula for cosine we receive from \eqref{eq:even03} that

\begin{subequations}\label{eq:even04}
\begin{eqnarray}
\sum_{{\ell}=1}^{k}\cos 2{\ell}\alpha  
&=& \frac12 \left( \cos\alpha\ \frac{\sin(2k+2)\alpha}{\sin\alpha} - \cos (2k+2)\alpha -1 \right) \label{eq:even04a} \\ 
&=& \frac12 \left( \frac{\cos\alpha\sin(2k+2)\alpha-\sin\alpha\cos (2k+2)\alpha }{\sin\alpha} -1 \right) \label{eq:even04b}\\ 
&=& \frac12 \left( \frac{\sin\left((2k+2)-1\right)\alpha}{\sin\alpha}  -1 \right)\\
&=& \frac12 \left( \frac{\sin(2k+1)\alpha}{\sin\alpha}  -1 \right). \label{eq:even04c}
\end{eqnarray}
\end{subequations}

If $\varphi=2\alpha$ then 
\begin{eqnarray}\label{eq:05}
\sum_{{\ell}=1}^{k}  \cos {\ell}\varphi =
  \frac12\left(  \frac{\sin(k+\frac12)\varphi}{\sin\frac12 \varphi}  -1\right).
\end{eqnarray}

\section{Other summation form}

In this section, we give an other summation form for the cosines without half angles by the help of the defined geometrical construction.

Now let us take the orthogonal projection of the  segments $A_{{\ell}-1}A_{\ell}$ $({\ell}=1,2,...,n)$ onto the line $e$ (see Figure \ref{fig:even} and \ref{fig:even_n}) and summarize them for all ${\ell}$ from 1 to $2k+1$. (The sum is equal to $x(\alpha)$ if $n=2k+1$ and $A_1\in e$.) Now we gain a similar equation to \eqref{eq:even03}, namely 

\begin{eqnarray*}\label{eq:odd01}
2\left( \cos \alpha + \cos 3\alpha+ \cdots +\cos (2k-1)\alpha \right)+ \cos (2k+1)\alpha = \cos\alpha\ \frac{\sin(2k+1)\alpha}{\sin\alpha}. 
\end{eqnarray*}
With analogous calculation to \eqref{eq:even04} we obtain
\begin{subequations}
\begin{eqnarray}
\sum_{{\ell}=1}^{k}\cos (2{\ell}-1)\alpha 
&=& \frac12 \left( \cos\alpha\ \frac{\sin(2k+1)\alpha}{\sin\alpha} - \cos (2k+1)\alpha  \right)\label{eq:odd02a}\\ 
&=& \frac12 \left( \frac{\sin 2k\alpha}{\sin\alpha}   \right)\label{eq:odd02b}.
\end{eqnarray}
\end{subequations}

Summing equations \eqref{eq:even04c} and \eqref{eq:odd02b}, we have
\begin{eqnarray*}
\sum_{{\ell}=1}^{k}\cos 2{\ell}\alpha + \sum_{{\ell}=1}^{k}\cos (2{\ell}-1)\alpha 
&=& \frac12 \left( \frac{\sin (2k+1)\alpha}{\sin\alpha}+ \frac{\sin 2k\alpha}{\sin\alpha} -1 \right),
\end{eqnarray*}  
and finally if $m=2k$ and $\varphi=\alpha$ we obtain formula \eqref{eq:sum2}.


\end{document}